\documentclass[twoside,a4paper,10pt]{article}
\usepackage{latexsym,amssymb}
\usepackage{mltex}
\usepackage{graphicx}
\usepackage[frenchb,english]{babel}
\usepackage[matrix,arrow]{xy}
\def \be{\begin{eqnarray*}}
\def \ee{\end{eqnarray*}}
\def \ben{\begin{enumerate}}
\def \een{\end{enumerate}}
\def \beit{\begin{itemize}}
\def \eeit{\end{itemize}}
\def \bui#1#2{\mathrel{\mathop{\kern 0pt#1}\limits^{#2}}}
\def \buil#1#2{\mathrel{\mathop{\kern 0pt#1}\limits_{#2}}}
\def \bfll{\begin{flushleft}}
\def \efll{\end{flushleft}}
\def \bflr{\begin{flushright}}
\def \eflr{\end{flushright}}

\def \findemo{\bflr$\Box$\eflr}
\def \lra{\longrightarrow}
\def \lmt{\longmapsto}
\def \ovl{\overline}
\def \wih{\widehat}
\def \wit{\widetilde}
\def \wnabla{\wit{\nabla}}
\def \cdotM{\buil{\cdot}{M}}

\newcommand{\pa}[1]{\left(#1\right)}
\newcommand{\la}{\langle}       
\newcommand{\ra}{\rangle} 

\newtheorem{edefi}{Def\mbox{}inition}
\newtheorem{elemme}{Lemma}
\newtheorem{erem}{Note}
\newtheorem{erems}{Notes}
\newtheorem{ecor}{Corollary}
\newtheorem{ethm}{Theorem}

\title{Dirac operators on Lagrangian submanifolds}
\author{Nicolas Ginoux}

\begin{document}
\maketitle

\noindent\begin{center}\begin{tabular}{p{155mm}}
\begin{small}{\bf Abstract.}
We study a natural Dirac operator on a Lagrangian submanifold of a K\"ahler manifold. We f\mbox{}irst show that its square coincides with the Hodge\,-\,de Rham Laplacian provided the complex structure identif\mbox{}ies the Spin structures of the tangent and normal bundles of the submanifold. We then give extrinsic estimates for the eigenvalues of that operator and discuss some examples.\end{small}\\
\end{tabular}\end{center}

\noindent\begin{small}{\it Mathematics Subject Classif\mbox{}ication}: 53C15, 53C27, 53C40 \\
\end{small}

\noindent\begin{small}{\it Key words}: Dirac operators, Global Analysis, Spectral Geometry, Spin Geometry, Lagrangian submanifolds\\
\end{small}

$ $\\

\noindent{\bf\Large Introduction\\}

The main object of this paper is to initiate the study of the properties of a Dirac operator on Lagrangian submanifolds of K\"ahler manifolds.\\

\noindent Spin geometry has revealed as a powerful tool in intrinsic geometry for a long time (see e.g. \cite{LM}). It is however a recent and striking fact that spinors play a role in extrinsic geometry as well. Initiated by Witten \cite{Wi81}, the use of Dirac operators on submanifolds has only been developed over the last years, especially about the following question: how can one relate analytical properties of some Dirac operators on a submanifold with extrinsic geometric quantities? For submanifolds of real space-forms, on which there exists particular spinor f\mbox{}ields (mainly parallel spinor f\mbox{}ields, up to a conformal change of the metric), a beautiful series of results has already appeared (see \cite{Ginthese} for references). However, answering that question in presence of further geometric structures seems to have been little considered.\\

\noindent We propose in this paper to begin with the study of (immersed) submanifolds of \emph{K\"ahlerian} manifolds. The presence of a complex structure on the ambient manifold gives rise to a rich variety of submanifolds (totally real, K\"ahlerian, real hypersurface,...). That is why we shall restrict our attention to a particular class of submanifolds, namely \emph{Lagrangian} submanifolds. A submanifold of a K\"ahlerian manifold is Lagrangian if and only if the (ambient) complex structure maps its tangent bundle \emph{onto} its normal bundle. Like every submanifold, a Lagrangian submanifold carries a \emph{twisted-Dirac operator}. We shall f\mbox{}irst prove that, if furthermore the complex structure identif\mbox{}ies the Spin structures of the tangent and normal bundles, then this twisted-Dirac operator identif\mbox{}ies with the \emph{Euler} operator. This requires adapting some technical algebraic Lemmas (compare with \cite{GinMor2002,Baer98}), which we shall therefore recall in detail in the f\mbox{}irst part. Coming back to the original question, we then prove new eigenvalue estimates for the above twisted-Dirac operator, and show their sharpness through examples. The results obtained show analogies with \cite{Chen83,Ros83}.\\

\noindent This work is partially based on the Ph.D.-thesis of the author \cite{Ginthese}.\\

\noindent {\bf Acknowledgement} It is a pleasure to thank Stefano Marchiafava for bringing the author's attention to the subject, Jean-Louis Milhorat and Sebasti\'an Montiel for fruitful discussions about it, Christian B\"ar for valuable comments, and Oussama Hijazi for his encouragements. The author would also like to thank the following institutions for their support: the Institut \'Elie Cartan de Nancy, the Max-Planck Institut f\"ur Mathematik in den Naturwissenschaften (Leipzig) and the University of Hamburg.\\

\section{Spin structures and Dirac operators on a Lagrangian submanifold}
We begin with collecting basic facts about Spin structures on Lagrangian submanifolds of K\"ahlerian manifolds (see also \cite{LM,BHMM,Friedlivre} for general Spin geometry). We f\mbox{}irst describe the necessary algebraic material, then transport it to bundles with the help of a group-equivariance condition.\\

\subsection{Clif\mbox{}ford algebras and spinors}
\noindent In this subsection, we recall some important isomorphisms between the complex Clif\mbox{}ford algebra (see def\mbox{}inition below) and other vector spaces. We point out that the isomorphism (\ref{IS2nSS}) below slightly dif\mbox{}fers from the equivalent one in \cite{GinMor2002} or \cite{Baer98}, since we want here to keep track of the ``Clif\mbox{}ford action'' in a more suitable way for our setting.\\
We f\mbox{}ix a positive integer $n$, and denote by ``$\mathrm{can}$'' the standard Euclidean inner product of $\mathbb{R}^n$. Throughout this paper, unless explicitely mentioned, all the isomorphisms will be denoted by the identity map.\\
Let $\mathbb{C}\mathrm{l}_n$ (resp. $\mathrm{Cl}_n$) be the complex (resp. real) Clif\mbox{}ford algebra of $(\mathbb{R}^n,\mathrm{can})$, that is, the only associative complex (resp. real) algebra with unit generated by $\mathbb{R}^n$ with the relation
\[v\cdot w\cdot+\;w\cdot v\cdot=-2\mathrm{can}(v,w)1,\] for all vectors $v$ and $w$ in $\mathbb{R}^n$. The product ``$\,\cdot\,$'' is called the \emph{Clif\mbox{}ford multiplication}. We recall the properties of $\mathbb{C}\mathrm{l}_n$ which will be important for the future:
\beit\item  Let $\Lambda\mathbb{R}^n\otimes\mathbb{C}$ be the complexif\mbox{}ied exterior algebra of $\mathbb{R}^n$. Then there exists a \emph{canonical} linear isomorphism \cite{LM}
\begin{equation}\label{IClifL}\mathbb{C}\mathrm{l}_n\lra\Lambda\mathbb{R}^n\otimes\mathbb{C}\end{equation} which maps every element of the form $v\cdot\varphi$ ($v\in\mathbb{R}^n$, $\varphi\in\mathbb{C}\mathrm{l}_n$) onto $v\wedge\varphi-v\lrcorner\,\varphi$, where ``$v\lrcorner\,\varphi$'' stands for $v^\flat\lrcorner\,\varphi$ through the musical isomorphism $v\lmt v^\flat:=\mathrm{can}(v,\cdot)$ between $\mathbb{R}^n$ and $(\mathbb{R}^n)^*$. We hereby identify through (\ref{IClifL}) the space $\Lambda^p\mathbb{R}^n\otimes\mathbb{C}$ as a subspace of $\mathbb{C}\mathrm{l}_n$.
\item The algebra $\mathbb{C}\mathrm{l}_n$ is either a matrix algebra or the copy of two such ones: there exists a complex vector space $\Sigma_n$ of dimension $2^{[\frac{n}{2}]}$, called the \emph{space of spinors}, and an isomorphism of complex algebras
\begin{equation}\label{ClifS}
\mathbb{C}\mathrm{l}_n\cong\left\{\begin{array}{ll}\mathrm{End}_{\mathbb{C}}(\Sigma_n)& \textrm{ if }n\textrm{ is even}\\
\mathrm{End}_{\mathbb{C}}(\Sigma_n)\oplus\mathrm{End}_{\mathbb{C}}(\Sigma_n)& \textrm{ if }n\textrm{ is odd.}\end{array}\right.
\end{equation}

\noindent Without loss of generality (see e.g. \cite{LM}), we further assume that, when $n$ is odd, the isomorphism (\ref{ClifS}) maps the complex volume element of $\mathbb{C}\mathrm{l}_n$ (see e.g. in \cite{BHMM} for its def\mbox{}inition) onto $\mathrm{Id}_{\Sigma_n}\oplus -\mathrm{Id}_{\Sigma_n}$. We def\mbox{}ine $\delta_n$ as the isomorphism (\ref{ClifS}) if $n$ is even, and the composition of the projection onto the f\mbox{}irst subalgebra $\mathrm{End}_{\mathbb{C}}(\Sigma_n)$ with (\ref{ClifS}) if $n$ is odd. In particular, for $n$ odd and for every $v$ in $\mathbb{R}^n$, the isomorphism (\ref{ClifS}) reads $v\lmt\delta_n(v)\oplus -\delta_n(v)$, see \cite{LM}.
\item The space $\Sigma_n$ carries a natural Hermitian inner product ``$\la\cdot\,,\cdot\ra$'' (which we assume to be complex-linear in the \emph{f\mbox{}irst} argument) such that, for every $v$ in $\mathbb{R}^n$ and all $\sigma,\sigma'$ in $\Sigma_n$,
\begin{equation}\label{MUantiH}\la\delta_n(v)\sigma,\sigma'\ra=-\la\sigma,\delta_n(v)\sigma'\ra.\end{equation}  
The property (\ref{MUantiH}) determines the Hermitian inner product ``$\la\cdot\,,\cdot\ra$'' up to a positive scalar \cite{BHMM}.\\

\eeit
\noindent Def\mbox{}ine now the \emph{Spin group} $\mathrm{Spin}_n$ as
\[\mathrm{Spin}_n:=\left\{v_1\cdot\ldots\cdot v_{2k}\quad /\quad k\geq 1, \quad v_j\in\mathbb{R}^n,\quad \mathrm{can}(v_j,v_j)=1\right\}\] and the \emph{Spin representation} to be the restriction of $\delta_n$ to $\mathrm{Spin}_n$. The Spin group is a compact Lie-subgroup of the group of invertible elements in $\mathbb{C}\mathrm{l}_n$ which has the following remarkable properties:
\beit\item There exists a two-fold covering Lie-group-homomorphism from $\mathrm{Spin}_n$ onto the special orthogonal group $\mathrm{SO}_n$, which we denote by ``$\mathrm{Ad}$''.
\item Denoting also ``$\mathrm{Ad}$'' the composition of the natural representation of $\mathrm{SO}_n$ on $\Lambda\mathbb{R}^n\otimes\mathbb{C}$ with $\mathrm{Ad}$, the isomorphism (\ref{IClifL}) is $\mathrm{Spin}_n$-equivariant, i.e., for every $u$ in $\mathrm{Spin}_n$ and $\varphi$ in $\mathbb{C}\mathrm{l}_n$, 
\[u\cdot\varphi\cdot u^{-1}\simeq \mathrm{Ad}(u)\varphi\] through (\ref{IClifL}).
\item Every Hermitian inner product ``$\la\cdot\,,\cdot\ra$'' satisfying (\ref{MUantiH}) is $\mathrm{Spin}_n$-invariant, i.e., the Spin representation is unitary w.r.t. ``$\la\cdot\,,\cdot\ra$''.\\\eeit

\noindent We now recall two lemmas and discuss their consequences.\\

\begin{elemme}\label{ljmath}
There exists a \emph{complex-antilinear} automorphism $\jmath$ of $\Sigma_n$ commuting with the Spin representation, i.e., for every $u$ in $\mathrm{Spin}_n$, we have $\delta_n(u)\circ\jmath=\jmath\circ\delta_n(u)$.\\
\end{elemme}

\noindent{\it Proof}: Although it follows from representation theory (see e.g. \cite{Tits67}, p. 21), we give here an elementary argument.\\
From the classif\mbox{}ication of \emph{real} Clif\mbox{}ford algebras (see \cite{LM}), we have :
\begin{displaymath}\mathrm{Cl}_n\cong\left\{\begin{array}{lll}\mathbb{R}\pa{2^{[\frac{n}{2}]}} & \textrm{ if }n\equiv 0\textrm{ or }6& (8)\\
\mathbb{R}\pa{2^{[\frac{n}{2}]}}\oplus \mathbb{R}\pa{2^{[\frac{n}{2}]}} & \textrm{ if }n\equiv 7& (8)\\
\mathbb{H}\pa{2^{[\frac{n-1}{2}]}} & \textrm{ if }n\equiv 2\textrm{ or }4& (8)\\
\mathbb{H}\pa{2^{[\frac{n-2}{2}]}}\oplus\mathbb{H}\pa{2^{[\frac{n-2}{2}]}} & \textrm{ if }n\equiv 3& (8)\\
\mathbb{C}\pa{2^{[\frac{n}{2}]}}& \textrm{ if }n\equiv 1& (4)\end{array}\right.\end{displaymath}

\noindent As $\mathbb{C}\mathrm{l}_n\cong\mathrm{Cl}_n\otimes\mathbb{C}$, we see that, if $n\equiv 6$, $7$ or $8 \;(8)$, the \emph{complex} representation $\delta_n:\mathbb{C}\mathrm{l}_n\lra \mathrm{End}_\mathbb{C}\pa{\Sigma_n}$ admits a real structure, i.e., there exists a $\mathbb{C}$-antilinear and involutive automorphism  $\jmath$ of $\Sigma_n$ such that, for every vector $v$ in $\mathbb{R}^n$ (hence for every element in $\mathbb{C}\mathrm{l}_n$), 
\[\delta_n(v)\circ\jmath=\jmath\circ\delta_n(v).\] If $n\equiv 2$, $3$ or $4 \;(8)$, there exists a quaternionic structure on $\Sigma_n$, i.e. a $\mathbb{C}$-antilinear automorphism $\jmath$ of $\Sigma_n$ satisfying $\jmath^2=-\mathrm{Id}$ and the preceding relation. If $n\equiv 1\;(4)$, the real representation of $\mathrm{Cl}_n$ being already complex, there exists no $\mathbb{C}$-antilinear automorphism  of $\Sigma_n$ commuting with the action of \emph{every} vector of $\mathbb{R}^n$ as before. However, that relation needs only to hold on $\mathrm{Spin}_n$ and not on $\mathrm{Cl}_n$. We eliminate the case $n=1$, for which this is obviously true. For $n>1$, as $\mathrm{Spin}_n$ is a subset of $\mathrm{Cl}_n^0:=\buil{\oplus}{p\textrm{ even}}\Lambda^p\mathbb{R}^n$, and $\mathrm{Cl}_n^0$ identif\mbox{}ies with $\mathrm{Cl}_{n-1}$ through an isomorphism which provides the equivalence of $\delta_{n-1}$ (or ``double copy'' as in (\ref{ClifS})) with $(\delta_n)_{|_{\mathbb{C}\mathrm{l}_n^0}}$ (see \cite{LM}), we just need to solve the problem for $\delta_{n-1}$. But from the preceding arguments, the representation $\delta_{n-1}$ admits such a structure, so that we again obtain a $\mathbb{C}$-antilinear automorphism $\jmath$ of $\Sigma_n$ which commutes with $(\delta_n)_{|_{\mathrm{Spin}_n}}$.\\

\noindent To sum up, for every $n\geq 1$, there exists a $\mathbb{C}$-antilinear automorphism $\jmath$ of $\Sigma_n$ such that, for every $u$ in $\mathrm{Spin}_n$,
\[\delta_n(u)\circ\jmath=\jmath\circ \delta_n(u),\] which is the desired property.\findemo

\begin{ecor}[\cite{Gaud2002}; \cite{Baer97}, p. 244]\label{cClifSS} There exists a complex-linear isomorphism
\begin{equation}\label{ClifSS}
\mathbb{C}\mathrm{l}_n\lra\left\{\begin{array}{ll}\Sigma_n\otimes\Sigma_n&\textrm{ if }n\textrm{ is even}\\
\Sigma_n\otimes\Sigma_n\oplus\Sigma_n\otimes\Sigma_n&\textrm{ if }n\textrm{ is odd}\end{array}\right.\end{equation} satisfying:
\beit\item For every $v$ in $\mathbb{R}^n$ and every $\varphi$ in $\mathbb{C}\mathrm{l}_n$, the element $v\cdot\varphi$ is mapped onto $\left\{\delta_n(v)\otimes\mathrm{Id}\right\}\varphi$ when $n$ is even (resp. onto $\left\{\delta_n(v)\otimes\mathrm{Id}\oplus -\delta_n(v)\otimes\mathrm{Id}\right\}\varphi$ when $n$ is odd),
\item The isomorphism (\ref{ClifSS}) is $\mathrm{Spin}_n$-equivariant: for every $u$ in $\mathrm{Spin}_n$ and every $\varphi$ in $\mathbb{C}\mathrm{l}_n$, the element $u\cdot\varphi\cdot u^{-1}$ is mapped onto $\{\delta_n(u)\otimes\delta_n(u)\}\varphi$ when $n$ is even (resp. onto $\{\delta_n(u)\otimes\delta_n(u)\oplus\,\delta_n(u)\otimes\delta_n(u)\}\varphi$ when $n$ is odd).\eeit\end{ecor}

\noindent{\it Proof}: The existence of $\jmath$ allows to def\mbox{}ine a \emph{complex}-linear isomorphism
\begin{eqnarray}\label{ISS*}
\Sigma_n&\lra&\Sigma_n^*\\
\nonumber \sigma&\lmt&\la\cdot\,,\jmath(\sigma)\ra,\end{eqnarray}
 which is $\mathrm{Spin}_n$-equivariant: for every $u$ in $\mathrm{Spin}_n$ and every $\sigma$ in $\Sigma_n$,
\be
\la\cdot\,,\jmath(\delta_n(u)\sigma)\ra&=&\la\cdot\,,\delta_n(u)\jmath(\sigma)\ra\\
&=&\la\delta_n(u)\delta_n(u^{-1})\cdot\,,\delta_n(u)\jmath(\sigma)\ra\\
&\bui{=}{(\ref{MUantiH})}&\la\delta_n(u^{-1})\cdot\,,\jmath(\sigma)\ra\\
&=&{}^t\pa{\delta_n(u^{-1})}\{\la\cdot\,,\jmath(\sigma)\ra\},\ee that is, the isomorphism (\ref{ISS*}) maps $\delta_n(u)\sigma$ onto ${}^t\pa{\delta_n(u^{-1})}\{\la\cdot\,,\jmath(\sigma)\ra\}$. Composing then the canonical isomorphism $\mathrm{End}_\mathbb{C}\pa{\Sigma_n}\cong\Sigma_n\otimes\Sigma_n^*$ with (\ref{ClifS}) and using (\ref{ISS*}), we obtain the isomorphism (\ref{ClifSS}), which obviously satisf\mbox{}ies the desired property w.r.t. the left Clif\mbox{}ford action. Let now $v$ in $\mathbb{R}^n$ such that $\mathrm{can}(v,v)=1$ and $\varphi$ in $\mathbb{C}\mathrm{l}_n$. Assume $n$ even and f\mbox{}ix an o.n.b. $(\sigma_k)_{1\leq k\leq n}$ of $\Sigma_n$; then, noting that (\ref{ISS*}) maps $\jmath^{-1}(\sigma_k)$ onto $\sigma_k^*$,
\be 
\sum_{k=1}^{2^{[\frac{n}{2}]}}\delta_n(\varphi\cdot v)\sigma_k\otimes\jmath^{-1}(\sigma_k)&=&\sum_{k=1}^{2^{[\frac{n}{2}]}}\delta_n(\varphi)\delta_n(v)\sigma_k\otimes\jmath^{-1}(\delta_n(-v)\delta_n(v)\sigma_k)\\
&=&-\sum_{k=1}^{2^{[\frac{n}{2}]}}\delta_n(\varphi)\delta_n(v)\sigma_k\otimes\delta_n(v)\jmath^{-1}(\delta_n(v)\sigma_k).
\ee
But, from (\ref{MUantiH}), for every unit vector $v$ of $\mathbb{R}^n$, the endomorphism $\delta_n(v)$ of $\Sigma_n$ is unitary, so that the vectors $\sigma_k':=\delta_n(v)\sigma_k$ form an o.n.b. of $\Sigma_n$. Hence we have
\be
\sum_{k=1}^{2^{[\frac{n}{2}]}}\delta_n(\varphi\cdot v)\sigma_k\otimes\jmath^{-1}(\sigma_k)&=&-\sum_{k=1}^{2^{[\frac{n}{2}]}}\delta_n(\varphi)\sigma_k'\otimes\delta_n(v)\jmath^{-1}(\sigma_k')\\
&=&-\{\mathrm{Id}\otimes\delta_n(v)\}\sum_{k=1}^{2^{[\frac{n}{2}]}}\delta_n(\varphi)\sigma_k'\otimes\jmath^{-1}(\sigma_k')\\
&=&-\{\mathrm{Id}\otimes\delta_n(v)\}\sum_{k=1}^{2^{[\frac{n}{2}]}}\delta_n(\varphi)\sigma_k\otimes\jmath^{-1}(\sigma_k),
\ee
for the second factor of the right member does not depend on the choice of the o.n.b. $(\sigma_k)_k$ of $\Sigma_n$. We hence have, through (\ref{ClifSS}),
\[\varphi\cdot v=-\{\mathrm{Id}\otimes\delta_n(v)\}\varphi,\] which then obviously holds for every $v$ in $\mathbb{R}^n$. For the case $n$ odd, the computations go the same way (just beware of the sign!) and we obtain $\varphi\cdot v=-\{\mathrm{Id}\otimes\delta_n(v)\oplus-\mathrm{Id}\otimes\delta_n(v)\}\varphi$, so that, for every $v$ in $\mathbb{R}^n$ and every $\varphi$ in $\mathbb{C}\mathrm{l}_n$,
\begin{equation}\label{MUdroite}
\varphi\cdot v\bui{=}{(\ref{ClifSS})}\left\{\begin{array}{ll}-\{\mathrm{Id}\otimes\delta_n(v)\}\varphi&\textrm{ if }n\textrm{ is even}\\
-\{\mathrm{Id}\otimes\delta_n(v)\oplus-\mathrm{Id}\otimes\delta_n(v)\}\varphi&\textrm{ if }n\textrm{ is odd.}\end{array}\right.\end{equation}
 Since, for a (real) unit vector $v$, we have $v^{-1}=-v$ in $\mathbb{C}\mathrm{l}_n\setminus\{0\}$, we deduce from (\ref{MUdroite}) that, for every $u$ in $\mathrm{Spin}_n$,
\[\varphi\cdot u\bui{=}{(\ref{ClifSS})}\left\{\begin{array}{ll}\{\mathrm{Id}\otimes\delta_n(u^{-1})\}\varphi&\textrm{ if }n\textrm{ is even}\\
\{\mathrm{Id}\otimes\delta_n(u^{-1})\oplus\mathrm{Id}\otimes\delta_n(u^{-1})\}\varphi&\textrm{ if }n\textrm{ is odd,}\end{array}\right.\]from which follows the $\mathrm{Spin}_n$-equivariance.\findemo

\noindent For the next lemma, we recall an explicit description of the space of spinors as a subspace of the complex Clif\mbox{}ford algebra in even dimensions (compare with \cite{BHMM,Kirch90,Kirch93}). We consider $\mathbb{R}^{2n}$ endowed with its natural complex structure $J$, so that  $\mathbb{R}^{2n}=\mathbb{R}^n\oplus J(\mathbb{R}^n)$. Let $p_\pm$ the two projectors of $\mathbb{R}^{2n}\otimes\mathbb{C}$ def\mbox{}ined by 
\[p_\pm:=\frac{1}{2}\pa{\mathrm{Id}\mp iJ},\] where $J$ is extended as a complex-linear automorphism of $\mathbb{R}^{2n}\otimes\mathbb{C}$. The endomorphisms $p_+$ and $p_-$ satisfy the following: $p_-\circ p_+=p_+\circ p_-=0$, $p_\pm\circ J=J\circ p_\pm=\pm ip_\pm$, and $\mathrm{can}(p_+(Z),Z')=\mathrm{can}(Z,p_-(Z'))$ for all $Z$, $Z'$ in $\mathbb{R}^{2n}\otimes\mathbb{C}$.\\
Let $(e_j)_{1\leq j\leq n}$ be the canonical basis of $\mathbb{R}^n$. For $1\leq j\leq n$, def\mbox{}ine $z_j:=p_+(e_j)$ and $\ovl{z}_j:=p_-(e_j)$. The vectors $z_1,\ldots,z_n,\ovl{z}_1,\ldots,\ovl{z}_n$ form the so-called \emph{Witt-basis} of $\mathbb{R}^{2n}\otimes\mathbb{C}$ associated to the basis $(e_j)_{1\leq j\leq n}$ of $\mathbb{R}^n$. Set
\[\ovl{\omega}:=\ovl{z}_1\cdot\ldots\cdot\ovl{z}_n.\]
From the above properties of $p_\pm$, that element of $\mathbb{C}\mathrm{l}_{2n}$ is independent of the choice of the p.o.n.b. of $\mathbb{R}^n$: replacing $(e_j)_{1\leq j\leq n}$ by another p.o.n.b. of $\mathbb{R}^n$, and taking the associated Witt-basis of $\mathbb{R}^{2n}\otimes\mathbb{C}$, one obtains the same element $\ovl{\omega}$.\\
For $1\leq p\leq n$, set $\mathcal{L}^p:=\mathrm{Span}_{\mathbb{C}}\left\{z_{i_1}\cdot\ldots\cdot z_{i_p}\cdot\ovl{\omega},\qquad 1\leq i_1<\cdots <i_p\leq n\right\}$ and $\mathcal{L}^0:=\mathbb{C}\ovl{\omega}$. Those subspaces of $\mathbb{C}\mathrm{l}_{2n}$ do not depend on the choice of a p.o.n.b. of $\mathbb{R}^n$, in the preceding sense. It can furthermore be shown that $\oplus_{p=0}^n\mathcal{L}^p$ is a left-ideal of dimension $2^n$ in $\mathbb{C}\mathrm{l}_{2n}$, hence is isomorphic to $\Sigma_{2n}$ (see \cite{BHMM}). We can then set
\[\Sigma_{2n}:=\oplus_{p=0}^n\mathcal{L}^p.\]
Through that identif\mbox{}ication, for each $\psi$ in $\mathbb{C}\mathrm{l}_{2n}$, the endomorphism $\delta_{2n}(\psi)$ is given by the left-Clif\mbox{}ford multiplication by $\psi$. For example (and this will be crucial for the future), the Clif\mbox{}ford multiplication by the K\"ahler form $\wit{\Omega}(\cdot\,,\cdot)=\mathrm{can}(J\cdot\,,\cdot)$ of $(\mathbb{R}^{2n},J)$ is given by 
\[\delta_{2n}(\wit{\Omega})=\oplus_{p=0}^ni(2p-n)\mathrm{Id}_{\mathcal{L}^p},\] that is, $\mathcal{L}^p$ is the eigenspace of $\delta_{2n}(\wit{\Omega})$ for the eigenvalue $i(2p-n)$.\\
Moreover, a Hermitian inner product satisfying (\ref{MUantiH}) can be def\mbox{}ined in the following way: for any $1\leq i_1<\ldots <i_p\leq n$ and $1\leq j_1<\ldots <j_q\leq n$, set
\[\la z_{i_1}\cdot\ldots\cdot z_{i_p}\cdot\ovl{\omega}\,,z_{j_1}\cdot\ldots\cdot z_{j_q}\cdot\ovl{\omega}\ra:=\left\{\begin{array}{ll}0& \textrm{ if }\{i_1,\ldots,i_p\}\neq\{j_1,\ldots,j_q\}\\
2^{[\frac{n+1}{2}]}& \textrm{ otherwise.}\end{array}\right.
\]
It can be shown that (\ref{MUantiH}) holds and that this Hermitian inner product does not depend on the choice of a p.o.n.b. of $\mathbb{R}^n$.\\

\noindent We furthermore def\mbox{}ine $\mathrm{Spin}_n'$ to be the Spin group of $(J\mathbb{R}^n,\mathrm{can})$, i.e.,
\[\mathrm{Spin}_n':=\left\{w_1\cdot\ldots\cdot w_{2k}\quad /\quad k\geq 1, \quad w_j\in J\mathbb{R}^n,\quad \mathrm{can}(w_j,w_j)=1\right\}\subset\mathbb{C}\mathrm{l}_{2n}.\] From the universal property of Clif\mbox{}ford algebras \cite{LM}, the linear isometry $J:\mathbb{R}^n\lra J\mathbb{R}^n$ induces a Lie-group-isomorphism $\wit{J}:\mathrm{Spin}_n\lra\mathrm{Spin}_n'$. We then set $\delta_n':=\delta_n{}_{|_{\mathrm{Spin}_n}}\circ (\wit{J})^{-1}:\mathrm{Spin}_n'\lra\mathrm{Aut}_\mathbb{C}(\Sigma_n)$.\\

\begin{elemme}\label{lS2n}
Consider $\mathbb{C}\mathrm{l}_n$ as canonically embedded in
$\mathbb{C}\mathrm{l}_{2n}$. Then the map
\begin{eqnarray}\label{IS2n}
\mathbb{C}\mathrm{l}_n&\lra&\Sigma_{2n}\\
\nonumber \varphi&\lmt&\varphi\cdot\ovl{\omega}
\end{eqnarray} is a complex-linear isomorphism satisfying:
\beit\item For every $v$ in $\mathbb{R}^n$ and every $\varphi$ in $\mathbb{C}\mathrm{l}_n$, the element $v\cdot\varphi$ is mapped onto $v\cdot\varphi\cdot\ovl{\omega}$,
 \item For every $w$ in $J\mathbb{R}^n$ and every $\varphi$ in $\Lambda^p\mathbb{R}^n\otimes\mathbb{C}$, the element $(-1)^{p+1}i\varphi\cdot J(w)$ is mapped onto $w\cdot\varphi\cdot\ovl{\omega}$.\eeit
In particular, the isomorphism (\ref{IS2n}) is $\mathrm{Spin}_n$-equivariant w.r.t. the ``diagonal immersion''
\begin{eqnarray}\label{immdiag}
\mathrm{Spin}_n&\lra&\mathrm{Spin}_{2n}\\
\nonumber u&\lmt& u\cdot\wit{J}(u),\end{eqnarray} that is, for every $u$ in $\mathrm{Spin}_n$ and $\varphi$ in $\mathbb{C}\mathrm{l}_n$, the element $u\cdot\varphi\cdot u^{-1}$ is mapped onto $u\cdot\wit{J}(u)\cdot\varphi\cdot\ovl{\omega}$.\\
\end{elemme}

\noindent{\it Proof}: Since the linear map (\ref{IS2n}) is obviously surjective (for each $1\leq p\leq n$ and each $1\leq i_1<\ldots <i_p\leq n$, the element $e_{i_1}\cdot\ldots\cdot e_{i_p}$ of $\mathbb{C}\mathrm{l}_n$ is mapped onto $z_{i_1}\cdot\ldots z_{i_p}\cdot\ovl{\omega}$), and both spaces have the same dimension, it is a linear isomorphism. Furthermore, (\ref{IS2n}) maps the subspace $\Lambda^p\mathbb{R}^n\otimes\mathbb{C}$ onto $\mathcal{L}^p$.\\
The f\mbox{}irst property is trivial. On the other hand, for every $w$ in $J\mathbb{R}^n$ and $\varphi$ in $\Lambda^p\mathbb{R}^n\otimes\mathbb{C}$,
\be
(-1)^{p+1}i\varphi\cdot J(w)\cdot\ovl{\omega}&=& (-1)^p\varphi\cdot (-ip_+(J(w))\cdot\ovl{\omega})\\
&=&(-1)^p\varphi\cdot p_+(w)\cdot\ovl{\omega}\\
&=&(-1)^p\varphi\cdot w\cdot\ovl{\omega}\\
&=&w\cdot\varphi\cdot\ovl{\omega}, 
\ee hence the second point holds. As a consequence, for all $w_1$, $w_2$ in $J\mathbb{R}^n$ with $\mathrm{can}(w_1,w_1)=\mathrm{can}(w_2,w_2)=1$ and every $\varphi$ in $\mathcal{L}^p$, the preimage through (\ref{IS2n}) of $w_1\cdot w_2\cdot\varphi\cdot\ovl{\omega}$ is given by 
\be
(-1)^pi\{w_2\cdot\varphi\}\cdot J(w_1)\cdot\ovl{\omega}&=&(-1)^p(-1)^{p+1}i^2\varphi\cdot J(w_2)\cdot J(w_1)\cdot\ovl{\omega}\\
&=&\varphi\cdot J(w_2)\cdot J(w_1)\cdot\ovl{\omega},\ee i.e., is equal to $\varphi\cdot J(w_2)\cdot J(w_1)$. Note that it does no longer depend on $p$. Since $J(w_2)\cdot J(w_1)=(\wit{J})^{-1}(w_2\cdot w_1)=(\wit{J})^{-1}\{(w_1\cdot w_2)^{-1}\}$, we obtain 
\[\varphi\cdot(\wit{J})^{-1}(u'^{-1})\bui{\simeq}{(\ref{IS2n})}u'\cdot\varphi\cdot\ovl{\omega}\] for every $u'$ in $\mathrm{Spin}_n'$ and every $\varphi$ in $\mathbb{C}\mathrm{l}_n$, from which follows the last statement.\findemo

\begin{ecor}\label{cS2nSS} There exists a complex-linear isomorphism
\begin{equation}\label{IS2nSS} \Sigma_{2n}\lra\left\{\begin{array}{ll}\Sigma_n\otimes\Sigma_n&\textrm{ if }n\textrm{ is even}\\
\Sigma_n\otimes\Sigma_n\oplus\Sigma_n\otimes\Sigma_n&\textrm{ if }n\textrm{ is odd,}\\\end{array}\right.\end{equation} satisfying:
\beit\item For every $v$ in $\mathbb{R}^n$ and $\varphi$ in $\Sigma_{2n}$, the element $\delta_{2n}(v)\varphi$ is mapped onto $\left\{\delta_n(v)\otimes\mathrm{Id}\right\}\varphi$ if $n$ is even (resp. onto $\left\{\delta_n(v)\otimes\mathrm{Id}\oplus -\delta_n(v)\otimes\mathrm{Id}\right\}\varphi$ if $n$ is odd),
\item For every $w$ in $J\mathbb{R}^n$ and $\varphi$ in $\mathcal{L}^p$, the element $\delta_{2n}(w)\varphi$ is mapped onto $(-1)^pi\left\{\mathrm{Id}\otimes\delta_n(J(w))\right\}\varphi$ if $n$ is even (resp. onto $(-1)^pi\left\{\mathrm{Id}\otimes\delta_n(J(w))\oplus -\mathrm{Id}\otimes\delta_n(J(w))\right\}\varphi$ if $n$ is odd),
\item If $\Sigma_n\otimes\Sigma_n$ is endowed with the tensor product of a Hermitian inner product satisfying (\ref{MUantiH}) with itself, then the isomorphism (\ref{IS2nSS}) is unitary.
\eeit
In particular, the inverse of (\ref{IS2nSS}) is $\mathrm{Spin}_n\times\mathrm{Spin}_n'$-equivariant w.r.t. the group-homomorphism
\be \mathrm{Spin}_n\times\mathrm{Spin}_n'&\lra&\mathrm{Spin}_{2n}\\
(u,u')&\lmt&u\cdot u',
\ee that is: for every $(u,u')$ in $\mathrm{Spin}_n\times\mathrm{Spin}_n'$ and every $\varphi$ in $\Sigma_{2n}$, the isomorphism (\ref{IS2nSS}) maps the element $\delta_{2n}(u\cdot u')\varphi$ onto $\left\{\delta_n(u)\otimes\delta_n'(u')\right\}\varphi$ if $n$ is even (resp. onto $\left\{\delta_n(u)\otimes\delta_n'(u')\oplus\delta_n(u)\otimes\delta_n'(u')\right\}\varphi$ if $n$ is odd).\end{ecor}

\noindent{\it Proof}: The isomorphism (\ref{IS2nSS}) is obtained bringing together the isomorphisms (\ref{ClifSS}) and (\ref{IS2n}), and straightforward satisf\mbox{}ies the f\mbox{}irst property. The second one follows from Lemma \ref{lS2n} and from (\ref{MUdroite}). The third one comes from the fact that the squared-norm of the image of $e_{i_1}\cdot\ldots\cdot e_{i_p}\cdot\ovl{\omega}$ is $2^{[\frac{n+1}{2}]}$ (remember that $(e_j)_{1\leq j\leq n}$ stands here for the canonical basis of $\mathbb{R}^n$). The last statement follows from the two f\mbox{}irst ones, since for all vectors $w_1$ and $w_2$ in $J\mathbb{R}^n$, the isomorphism (\ref{IS2nSS}) maps $\delta_{2n}(w_1)\delta_{2n}(w_2)$ onto $\mathrm{Id}\otimes\{\delta_n(J(w_1))\delta_n(J(w_2))\}$ if $n$ is even (resp. onto $\mathrm{Id}\otimes\{\delta_n(J(w_1))\delta_n(J(w_2))\}\oplus\mathrm{Id}\otimes\{\delta_n(J(w_1))\delta_n(J(w_2))\}$ if $n$ is odd).\findemo

\begin{ecor}\label{cS2nL} There exists a complex-linear isomorphism
\begin{equation}\label{IS2nL}
\Sigma_{2n}\lra\Lambda\mathbb{R}^n\otimes\mathbb{C}
\end{equation}
satisfying:
\beit\item For every $v$ in $\mathbb{R}^n$ and $\varphi$ in $\Sigma_{2n}$, the element $\delta_{2n}(v)\varphi$ is mapped onto $v\wedge\varphi-v\lrcorner\,\varphi$,
\item For every $w$ in $J\mathbb{R}^n$ and $\varphi$ in $\Sigma_{2n}$, the element $\delta_{2n}(w)\varphi$ is mapped onto $-i\{J(w)\wedge\varphi+J(w)\lrcorner\,\varphi\}$.
\eeit
In particular, the isomorphism (\ref{IS2nL}) is $\mathrm{Spin}_n$-equivariant w.r.t. (\ref{immdiag}), i.e., for every $u$ in $\mathrm{Spin}_n$ and every $\varphi$ in $\Sigma_{2n}$, the isomorphism (\ref{IS2nL}) maps $u\cdot\wit{J}(u)\cdot\varphi$ onto $\mathrm{Ad}(u)\varphi$.\end{ecor}

 \noindent{\it Proof}: As before, the isomorphism (\ref{IS2nL}) is obtained from the isomorphisms (\ref{IClifL}) and (\ref{IS2n}), and satisf\mbox{}ies the f\mbox{}irst property. The second statement comes from the fact that, for every vector $v$ in $\mathbb{R}^n$ and every $\varphi$ in $\Lambda^p\mathbb{R}^n\otimes\mathbb{C}$, the element $\varphi\cdot v$ corresponds through (\ref{IClifL}) to the form $(-1)^p\{v\wedge\varphi+v\lrcorner\,\varphi\}$ (see \cite{LM}). The last one straightforward follows from the preliminary remarks and Lemma \ref{lS2n}.\findemo

\subsection{Spinor bundles on a Lagrangian submanifold}

\noindent We now deduce from the f\mbox{}irst subsection isomorphisms between spinor bundles on a submanifold. Since we shall need those isomorphisms in the setting of Lagrangian submanifolds (see def\mbox{}inition below), we restrict to the case of a Riemannian submanifold $(M^n,g)$ of dimension $n$ immersed in a Riemannian manifold $(\wit{M}^{2n},g)$ of (real) dimension $2n$. We shall always use the following notations: $II$ will be the bundle-valued second fundamental form of the immersion, $H$ its mean-curvature vector f\mbox{}ield (in our convention, $H:=\frac{1}{n}\mathrm{tr}(II)$), and $\nabla^M$ (resp. $\wnabla$) will be the Levi-Civit\`a connection of $(M,g)$ (resp. of $(\wit{M},g)$). The induced covariant derivative on the exterior bundle will be denoted analogously.\\ We assume $\wit{M}$ to be Spin, and f\mbox{}ix a Spin structure $\mathrm{Spin}(T\wit{M})\lra\mathrm{SO}(T\wit{M})$. As it is in general impossible to induce a Spin structure from $\wit{M}$ to $M$ (compare with the case of an oriented hypersurface \cite{LM,Baer98,Morel2001}), we assume $M$ to be Spin as well and f\mbox{}ix a Spin structure $\mathrm{Spin}(TM)\lra\mathrm{SO}(TM)$ on $M$. Then the normal bundle $NM$ of $M$ in $\wit{M}$ is Spin, and carries an induced Spin structure, $\mathrm{Spin}(NM)\lra\mathrm{SO}(NM)$, for which there exists a principal-bundle-homomorphism $\mathrm{Spin}(TM)\times_M\mathrm{Spin}(NM)\lra\mathrm{Spin}(T\wit{M})_{|_M}$ making the following diagram commutative \cite{Mil65}:

$$\xymatrix{\mathrm{Spin}(TM)\times_M\mathrm{Spin}(NM)\ar[dd]\ar[r]&\mathrm{Spin}(T\wit{M})_{|_M}\ar[dd]\ar[dr]&\\  & &  M\\ \mathrm{SO}(TM)\times_M\mathrm{SO}(NM)\ar[r]&\mathrm{SO}(T\wit{M})_{|_M}\ar[ur]&}$$

\noindent Let  $\Sigma M$ (resp. $\Sigma N$, $\Sigma\wit{M}$) be the spinor bundle of $TM$ (resp. of $NM$, $T\wit{M}$), i.e., the complex vector bundle associated to the Spin bundle through the Spin representation. There are three fundamental objects on the spinor bundle:
\beit\item The isomorphism (\ref{ClifS}) being obviously $\mathrm{Spin}_n$-equivariant induces a bilinear map, called the \emph{Clif\mbox{}ford multiplication} 
\be
TM\times_M\Sigma M&\lra&\Sigma M\\
(X,\varphi)&\lmt&\gamma_M(X)\varphi\ee satisfying
\[\gamma_M(X)\gamma_M(Y)+\gamma_M(Y)\gamma_M(X)=-2g(X,Y)\mathrm{Id}_{\Sigma M}\] for all vectors $X$ and $Y$ in $TM$. The same holds for $\Sigma N$ and $\Sigma \wit{M}$; we denote by $X\cdot\varphi:=\gamma_{\wit{M}}(X)\varphi$ the Clif\mbox{}ford multiplication by a vector $X$ on an element $\varphi$ of $\Sigma\wit{M}$.
\item The spinor bundle $\Sigma M$ also inherits from the space of spinors a Hermitian inner product ``$\la\cdot\,,\cdot\ra_M$'' satisfying
\[\la\gamma_M(X)\varphi,\psi\ra_M=-\la\varphi,\gamma_M(X)\psi\ra_M\] for every $X$ in $TM$ and all $\varphi$, $\psi$ in $\Sigma M$. The same property holds for $\Sigma N$ and for $\Sigma\wit{M}$, for which such a Hermitian inner product will be denoted by ``$\la\cdot\,,\cdot\ra$''.
\item The Levi-Civit\`a connection of $(TM,g)$ induces a covariant derivative $\nabla^{\Sigma M}$ on $\Sigma M$ \cite{BHMM,LM}. This covariant derivative is metric w.r.t. ``$\la\cdot\,,\cdot\ra_M$'' and satisf\mbox{}ies the Leibniz rule w.r.t. the Clif\mbox{}ford multiplication. We denote by  $\nabla^{\Sigma N}$ (resp. $\wnabla$) that covariant derivative on $\Sigma N$ (resp. on $\Sigma\wit{M}$).\eeit 

\noindent We now compare the dif\mbox{}ferent spinor bundles on the submanifold $M$. We need further notations in that purpose. For a tangent vector $X$ to $M$ and an element $\phi$ of $\Sigma M\otimes\Sigma N$ if $n$ is even (resp. of $\Sigma M\otimes\Sigma N\oplus\Sigma M\otimes\Sigma N$ if $n$ is odd), we def\mbox{}ine ``$X\cdotM\phi$'' to be
\begin{displaymath}\left|\begin{array}{ll}\{\gamma_M(X)\otimes\mathrm{Id}_{\Sigma N}\}\phi& \textrm{if }n\textrm{ is even}\\ \{\gamma_M(X)\otimes\mathrm{Id}_{\Sigma N}\oplus -\gamma_M(X)\otimes\mathrm{Id}_{\Sigma N}\}\phi& \textrm{if }n\textrm{ is odd.}\end{array}\right.\end{displaymath} We furthermore set
\begin{displaymath}\nabla:=\left|\begin{array}{ll} \nabla^{\Sigma M\otimes\Sigma N} & \textrm{if }m\textrm{ or }n\textrm{ is even}\\ \nabla^{\Sigma M\otimes\Sigma N}\oplus \nabla^{\Sigma M\otimes\Sigma N} & \textrm{otherwise.}\end{array}\right.\end{displaymath} Note that $\nabla$ is \emph{not} the natural covariant derivative of $\Sigma M$, since from its def\mbox{}inition it depends on the covariant derivative of the normal bundle.\\
From the above homomorphism between Spin bundles and the preceding subsection, we have the following:

\begin{elemme}\label{lfS2nSS}
There exists a complex-vector bundle isomorphism
\begin{equation}\label{isomspin}\Sigma\wit{M}_{|_M}\lra\left\{\begin{array}{ll}\Sigma M\otimes\Sigma N& \textrm{if }n\textrm{ is even}\\\Sigma M\otimes\Sigma N\oplus\Sigma M\otimes\Sigma N& \textrm{if }n\textrm{ is odd}\end{array}\right.\end{equation}
satisfying:
\beit
\item For every tangent vector f\mbox{}ield $X$ on $M$ and every section $\phi$ of $\Sigma\wit{M}_{|_M}$, the isomorphism (\ref{isomspin}) maps the section $X\cdot\phi$ onto $X\cdotM\phi$,
\item  For every tangent vector f\mbox{}ield $X$ on $M$ and every section $\phi$ of $\Sigma\wit{M}_{|_M}$, 
\[\wnabla_X\phi=\nabla_X\phi+\frac{1}{2}\sum_{j=1}^ne_j\cdot II(X,e_j)\cdot\phi\] in any local o.n.b. $(e_j)_{1\leq j\leq n}$ of $TM$.
\eeit
Furthermore, the isomorphism (\ref{isomspin}) can be assumed to be unitary.\\
\end{elemme}

\noindent{\it Proof}: From its equivariance under the action of $\mathrm{Spin}_n\times\mathrm{Spin}_n'$, the isomorphism (\ref{IS2nSS}) straightforward induces the isomorphism (\ref{isomspin}) between the vector bundles. The f\mbox{}irst property is just the translation of that of (\ref{IS2nSS}) on vector bundles. The second one is deduced in a quite analogous way as in \cite{GinMor2002} from the three following points: use the local expressions of the covariant derivatives $\wnabla$ and $\nabla$ \cite{BHMM,LM}, apply the Gau\ss{}-Weingarten formula on $T\wit{M}_{|_M}$, and use the correspondence through (\ref{isomspin}) between the Clif\mbox{}ford multiplications by 2-forms, that is: for all vectors $X_1$ and $X_2$ in $TM$,
\[X_1\cdot X_2\cdot\bui{\simeq}{(\ref{isomspin})} X_1\cdotM X_2\cdotM,\] and for all vectors $\nu_1$, $\nu_2$ in $NM$,
\[\nu_1\cdot\nu_2\cdot\bui{\simeq}{(\ref{isomspin})}\mathrm{Id}\otimes\gamma_N(\nu_1)\gamma_N(\nu_2)\;(\oplus\,\mathrm{Id}\otimes\gamma_N(\nu_1)\gamma_N(\nu_2)),\] where the parentheses stand for the case ``$n$ odd'' (see previous subsection). The last remark is also a direct consequence of Corollary \ref{cS2nSS}.\findemo

\noindent In this setup, the most natural Dirac operators that can be introduced on the manifold $M$ are the so-called \emph{twisted-Dirac operator} $D_M^{\Sigma N}$ \cite{LM} and the \emph{Dirac-Witten operator} $\wih{D}$ \cite{Wi81}, respectively def\mbox{}ined in a local o.n.b $(e_j)_{1\leq j\leq m}$ of $TM$ by 
\[D_M^{\Sigma N}:=\sum_{j=1}^me_j\cdotM\nabla_{e_j},\qquad \wih{D}:=\sum_{j=1}^me_j\cdot\wnabla_{e_j}.\] Both operators, which act on the sections of $\Sigma\wit{M}_{|_M}$, are elliptic, and from Lemma \ref{lfS2nSS} are related by
\[\wih{D}=D_M^{\Sigma N}-\frac{mH}{2}\cdot.\] 
Furthermore, the operator $D_M^{\Sigma N}$ is formally self-adjoint (but $\wih{D}$ is \emph{not}).\\

\noindent We now specialize to submanifolds with particular geometric structures. It is f\mbox{}irst important to point out that the spinor bundle $\Sigma N$ is in general \emph{not} isomorphic to $\Sigma M$; this may hold even if there exists an isomorphism between $TM$ and $NM$, such as for Lagrangian submanifolds in K\"ahlerian manifolds (see the examples in Notes \ref{rspinDir}). We therefore recall the notion of isomorphism between Spin structures:

\begin{edefi}\label{disomSpin}
Let $E$ and $F$ be two Spin vector bundles on a manifold $M$, with f\mbox{}ixed Spin structures $\mathrm{Spin}(E)\lra\mathrm{SO}(E)$ and $\mathrm{Spin}(F)\lra\mathrm{SO}(F)$. An \emph{isomorphism} between the Spin structures of $E$ and $F$ is given by a pair of principal-bundle isomorphisms $\mathrm{Spin}(E)\bui{\lra}{\wit{f}}\mathrm{Spin}(F)$ and $\mathrm{SO}(E)\bui{\lra}{f}\mathrm{SO}(F)$ such that the following diagram commutes:
$$\xymatrix{\mathrm{Spin}(E)\ar[dd]\ar[r]^{\wit{f}}&\mathrm{Spin}(F)\ar[dd]\ar[dr]&\\  & &  M\\ \mathrm{SO}(E)\ar[r]^f&\mathrm{SO}(F)\ar[ur]&}$$
\end{edefi}

\noindent If two vector bundles have isomorphic Spin structures, they obviously have isomorphic spinor bundles as well. Hence we give the following

\begin{ecor}\label{cfS2nL}
Assume that there exists an orientation-preserving isometry $f$ from $TM$ to $NM$ which induces an isomorphism $(\wit{f},f)$ of the respective Spin structures. Then there exists a complex-vector bundle isomorphism
\begin{equation}\label{IfS2nL}\Sigma\wit{M}_{|_M}\lra\Lambda TM\otimes\mathbb{C}\end{equation} satisfying:
\beit\item For every tangent vector $X$ to $M$ and every $\phi$ in $\Sigma\wit{M}_{|_M}$, the element $X\cdot\phi$ is mapped onto $X\wedge\phi-X\lrcorner\,\phi$,
\item For every vector $\nu$ in $NM$ and every $\phi$ in $\Sigma\wit{M}_{|_M}$, the element $\nu\cdot\phi$ is mapped onto $i\{f^{-1}(\nu)\wedge\phi+f^{-1}(\nu)\lrcorner\,\phi\}$,
\item If furthermore $f$ is parallel w.r.t. the respective connections on $TM$ and $NM$, then for every tangent vector f\mbox{}ield $X$ to $M$ and every section $\phi$ of $\Sigma\wit{M}_{|_M}$, the element $\nabla_X\phi$ is mapped onto $\nabla_X^M\phi$.
\eeit
\end{ecor} 

\noindent{\it Proof}: The existence of (\ref{IfS2nL}) is a direct consequence of Corollary \ref{cS2nL} and the fact that the Spin structure of $T\wit{M}_{|_M}$ reduces via $f$ to the Spin structure of $TM$. The f\mbox{}irst property comes straigthforward. For the second one, it is to be noted that the automorphism
\[J:=\left(\begin{array}{ll}0&-f^{-1}\\
f&0\end{array}\right)\] of $T\wit{M}_{|_M}$ is described through the $\mathrm{Spin}_n$-reduction as
\be \mathrm{Spin}(TM)\times_{\mathrm{Ad}}\mathbb{R}^{2n}&\lra& \mathrm{Spin}(TM)\times_{\mathrm{Ad}}\mathbb{R}^{2n}\\
 \lbrack\wit{s},v\rbrack&\lmt&\lbrack\wit{s},J(v)\rbrack.\ee The last statement follows from a short computation using the properties of compatibility between $\nabla$ and the other objects on $\Sigma M\otimes\Sigma M\;(\oplus\Sigma M\otimes\Sigma M)$.\findemo

\noindent Remark that, from the preceding proof, the existence of an orientation-preserving isometry $f:TM\lra NM$ is equivalent to the existence of an almost-Hermitian structure $J$ on $T\wit{M}_{|_M}$ mapping $TM$ onto $NM$. Let $\wit{\Omega}$ then denote the K\"ahler form of $(T\wit{M}_{|_M},g,J)$, i.e., $\wit{\Omega}(X,Y):=g(J(X),Y)$ for all $X$ and $Y$ in $T\wit{M}_{|_M}$. Under the hypotheses of Corollary \ref{cfS2nL}, the following holds: for every $0\leq p\leq n$ and $\phi$ in $\Lambda^pTM\otimes\mathbb{C}$,
\[\wit{\Omega}\cdot\phi=i(2p-n)\phi\] through (\ref{IfS2nL}). This also follows from the properties of (\ref{IS2n}), see previous subsection.\\

\noindent The existence of an almost-complex structure on $T\wit{M}_{|_M}$ is precisely the case we shall be interested in, since we shall consider submanifolds of K\"ahlerian manifolds. We now recall the following

\begin{edefi}\label{dSousvL}
A submanifold $M^n$ of a K\"ahlerian manifold $(\wit{M}^{2n},g,J)$ is called \emph{Lagrangian} if and only if 
\[J(TM)=NM,\] i.e., the complex structure identif\mbox{}ies the tangent and normal bundles of the submanifold.
\end{edefi}

\noindent For a Lagrangian submanifold in a K\"ahlerian manifold, the complex structure $J$ obviously preserves the metric and the orientation of $T\wit{M}_{|_M}$, and is parallel.

\begin{ecor}\label{cfS2nLL}
Let $(M^n,g)$ be a Spin Lagrangian submanifold immersed in a K\"ahlerian Spin manifold $(\wit{M}^{2n},g,J)$. Let the normal bundle $NM$ carry the induced Spin structure. Assume that the complex structure $J$ induces an isomorphism between the Spin structures of $TM$ and $NM$. Then there exists a complex-vector bundle isomorphism
\[\Sigma\wit{M}_{|_M}\lra\Lambda TM\otimes\mathbb{C}\] satisfying:
\beit\item For every tangent vector $X$ to $M$ and every $\phi$ in $\Sigma\wit{M}_{|_M}$, the element $X\cdot\phi$ is mapped onto $X\wedge\phi-X\lrcorner\,\phi$,
\item For every vector $\nu$ in $NM$ and every $\phi$ in $\Sigma\wit{M}_{|_M}$, the element $\nu\cdot\phi$ is mapped onto $-i\{J(\nu)\wedge\phi+J(\nu)\lrcorner\,\phi\}$,
\item For every tangent vector f\mbox{}ield $X$ to $M$ and every section $\phi$ of $\Sigma\wit{M}_{|_M}$, the element $\nabla_X\phi$ is mapped onto $\nabla_X^M\phi$.
\eeit
In particular, for each $0\leq p\leq n$, the subspace $\Lambda^p TM\otimes\mathbb{C}$ is the eigenspace associated to the eigenvalue $i(2p-n)$ of the action of the K\"ahler form $\wit{\Omega}$ of $(\wit{M}^{2n},g,J)$. Furthermore, for every section $\phi$ of $\Sigma\wit{M}_{|_M}$,
\[D_M^{\Sigma N}\phi=(d+\delta)\phi\qquad\textrm{ and }\qquad  \wih{D}\phi=(d+\delta)\phi+\frac{im}{2}\{J(H)\wedge\phi+J(H)\lrcorner\,\phi\},\] where $d$ (resp. $\delta$) denotes the exterior dif\mbox{}ferential (resp. codif\mbox{}ferential).\\
\end{ecor}

\noindent{\it Proof}: The only statement to be proved is the last one, for which it suf\mbox{}f\mbox{}ices to know that, for any local o.n.b. $(e_j)_{1\leq j\leq n}$ of $TM$, 
\[d=\sum_{j=1}^ne_j\wedge\nabla_{e_j}^M\qquad\textrm{ and }\qquad\delta=-\sum_{j=1}^ne_j\lrcorner\,\nabla_{e_j}^M.\]\findemo

\noindent \begin{erems}\label{rspinDir}\begin{rm}$ $\\
\noindent 1) In the same way as above, one can give a ``bundle-version'' of Corollary \ref{cClifSS}: let $E$ be any arbitrary Riemannian Spin vector bundle on a Spin manifold $M$ such that there exists an isomorphism from $TM$ to $E$, preserving the metric, the orientation \emph{and the Spin structure}. Then there exists a complex-vector bundle isomorphism between the Clif\mbox{}ford bundle and the tensor product $\Sigma M\otimes\Sigma E$ (or double copy), mapping $X\cdot\phi$ onto $X\cdotM\phi$ for every $X$ in $TM$ and $\phi$ in the Clif\mbox{}ford bundle; if furthermore the isomorphism from $TM$ to $E$ is parallel w.r.t. the covariant derivatives on $TM$ and $E$, then $\nabla_X^M\phi$ is mapped onto $\nabla_X^{\Sigma M\otimes\Sigma E}\phi$ (or double copy).\\
\noindent 2) We proved in Corollary \ref{cfS2nLL} that, under a compatibility condition between the complex structure and the Spin structures of $TM$ and $NM$, the twisted-Dirac operator can be identif\mbox{}ied with $d+\delta$ (the so-called \emph{Euler operator}), that is, a square-root of the Hodge\,-\,de Rham Laplacian $d\delta+\delta d$. This compatibility hypothesis is important, since otherwise the conclusions of Corollary \ref{cfS2nLL} may fail as can be seen on the following example. Consider the unit circle $M:=S^1$, canonically embedded in the complex line $\wit{M}:=\mathbb{C}$. This embedding is isometric and Lagrangian. Furthermore, $S^1$ carries two Spin structures, a trivial one and a non-trivial one. If one chooses the trivial (resp. non-trivial) Spin structure on the tangent bundle of $S^1$, then the induced Spin structure on the normal bundle is non-trivial (resp. trivial) \cite{BaerBeyr2001,Friedlivre}. Therefore, the complex structure does not even preserve the Spin bundles over $S^1$. Furthermore, the induced twisted Dirac operator is in both cases the fundamental Dirac operator of $S^1$ for the non-trivial Spin structure. Since this operator has trivial kernel, it cannot coincide with a square-root of the Hodge\,-\,de Rham Laplacian. One therefore sees that the hypothesis of compatibility of Corollary \ref{cfS2nLL} between the complex and the Spin structures is necessary.\\
\end{rm}\end{erems}

\section{An upper eigenvalue bound for the twisted Dirac operator on a Lagrangian submanifold}

\noindent In this section, we consider a compact Lagrangian Spin submanifold $(M^n,g)$ in a K\"ahlerian Spin manifold $(\wit{M}^{2n},g,J)$. Since the operator $D_M^{\Sigma N}$ is elliptic and formally self-adjoint, it has a discrete spectrum; we then denote by $\lambda_k$ ($k\in\mathbb{N}\setminus\{0\}$) its eigenvalues, counted with their multiplicities, assuming that $|\lambda_{k+1}|\geq |\lambda_k|$ for every $k\geq 1$.\\
We are interested in the following question: how can one control the smallest eigenvalues of the twisted Dirac operator in terms of extrinsic geometric invariants? For submanifolds of certain real space-forms, it was proved by C. B\"ar in \cite{Baer98} and the author in \cite{Ginthese,GinHeinD2003} that the ambient curvature together with either the $L^2$ or the $L^\infty$ norm of the mean curvature appear as the best candidates in that purpose. Those results were obtained considering restrictions to the submanifold of particular spinor f\mbox{}ields on the ambient manifold, called \emph{Killing spinors} (see \cite{BFGK} about those). As non Ricci-f\mbox{}lat K\"ahlerian Spin manifolds of (real) dimension greater than 2 do not admit such spinor f\mbox{}ields \cite{Hij86,Lich86,Lich87}, it comes as a natural question whether such kind of estimates could still hold in our context. We give an af\mbox{}f\mbox{}irmative and sharp answer to that problem, using the notion of \emph{K\"ahlerian Killing 
spinors} introduced by K.-D. Kirchberg in \cite{Kirch93} and O. Hijazi in \cite{Hij94}. Remember that, for a complex constant $\alpha$, an \emph{$\alpha$-K\"ahlerian Killing spinor} on the K\"ahlerian Spin manifold $(\wit{M}^{2n},g,J)$ is a couple of sections $(\psi,\phi)$ of $\Sigma\wit{M}$ satisfying, for every tangent vector f\mbox{}ield $Z$ on $\wit{M}$,
\begin{displaymath}\left\{\begin{array}{lll}\wnabla_Z\psi+\alpha p_-(Z)\cdot\phi&=& 0\\ \wnabla_Z\phi+\alpha p_+(Z)\cdot\psi&=& 0,\end{array}\right.\end{displaymath} where $p_\pm(Z):=\frac{1}{2}\pa{Z\mp iJ(Z)}$. When $\alpha=0$, an $\alpha$-K\"ahlerian Killing spinor is just a pair of parallel spinor f\mbox{}ields. As for Killing spinors, the presence of non-zero K\"ahlerian Killing spinors yields strong conditions on the geometry of $\wit{M}$ (see \cite{Kirch93,Hij94}): if $\alpha\neq 0$, the complex dimension $n$ of $\wit{M}$ has to be odd, the manifold $(\wit{M},g,J)$ has to be Einstein with scalar curvature $n(n+1)\alpha^2$ (therefore $\alpha$ must be either real or purely imaginary), and the sections $\psi$ and $\phi$ have to lie in particular eigenspaces of the Clif\mbox{}ford action of the K\"ahler form $\wit{\Omega}$ of $(\wit{M},g,J)$:
\be \wit{\Omega}\cdot\psi&=& -i\psi\\
\wit{\Omega}\cdot\phi&=& i\phi\ee
\noindent (remember that, in our convention, $\wit{\Omega}(X,Y):=g\pa{J(X),Y}$ for all vectors $X$ and $Y$ in $T\wit{M}$).

\subsection{Main result}

\noindent From here on, we denote by $\mathcal{K}_\alpha$ the space of $\alpha$-K\"ahlerian Killing spinors on $(\wit{M},g,J)$ (note that, if $\alpha\neq 0$, then $\mathcal{K}_\alpha\cap\mathcal{K}_{-\alpha}=\{0\}$). Manifolds carrying a non-zero $\mathcal{K}_\alpha$ have been completely characterized by A. Moroianu in \cite{Moroi95} when $\alpha$ is a non-zero real number, and partially by K.-D. Kirchberg in \cite{Kirch93} and M. Herzlich in \cite{Herz98} when $\alpha$ is purely imaginary. We prove the following:\\

\begin{ethm}\label{tmajLagr}
Let $(M^n,g)$ be a Lagrangian Spin submanifold of a K\"ahlerian Spin manifold $(\wit{M}^{2n},g,J)$. Let the normal bundle of $M$ in $\wit{M}$ carry the induced Spin structure. Assume that, for a given complex constant $\alpha$, the dimension of $\mathcal{K}_\alpha$ is $N$ ($N\geq 1$). Then the N$^{\textrm{th}}$ eigenvalue $\lambda_N$ of the twisted-Dirac operator $D_M^{\Sigma N}$ satisf\mbox{}ies
\[\lambda_N^2\leq\frac{(n+1)^2\alpha^2}{4}+\left\{\begin{array}{ll}\frac{n^2}{4\mathrm{Vol}(M)}\int_M|H|^2v_g&\textrm{ if }\alpha\textrm{ is real}\\ & \\ \frac{n^2||H||_\infty^2}{4}&\textrm{ if }\alpha\textrm{ is purely imaginary,}\end{array}\right.\] where $H$ is the mean curvature vector f\mbox{}ield of $M$ in $\wit{M}$.
\end{ethm}

$ $\\

\noindent{\it Proof}: Let $(\psi,\phi)$ be a non-zero $\alpha$-K\"ahlerian Killing spinor as above on $\wit{M}$. We compute the Rayleigh-quotient
\[\mathcal{Q}\pa{(D_M^{\Sigma N})^2,\psi+\phi}:=\frac{\int_M\la(D_M^{\Sigma N})^2(\psi+\phi)\,,\psi+\phi\ra v_g}{\int_M\la\psi+\phi\,,\psi+\phi\ra v_g}\] and apply the Min-Max principle. To obtain $(D_M^{\Sigma N})^2(\psi+\phi)$, we f\mbox{}irst evaluate $\wih{D}^2$ on $\psi+\phi$, then use the following relation (\cite{Ginthese}, Lemme 4.1): for every section $\varphi$ of $\Sigma\wit{M}_{|_M}$ and in every local o.n.b. $(e_j)_{1\leq j\leq n}$ of $TM$,
\begin{equation}\label{esqwihDDt}
(D_M^{\Sigma N})^2\varphi=\wih{D}^2\varphi+\frac{n^2|H|^2}{4}\varphi+\frac{n}{2}\sum_{j=1}^ne_j\cdot \nabla_{e_j}^NH\cdot\varphi,
\end{equation}
where $\nabla^NH$ denotes the normal covariant derivative of $H$.\\

\noindent Let us f\mbox{}ix a local o.n.b. $(e_j)_{1\leq j\leq n}$ of $TM$. From the hypotheses,
\be
\wih{D}\psi&=& \sum_{j=1}^ne_j\cdot\wnabla_{e_j}\psi\\
&=&-\alpha\sum_{j=1}^ne_j\cdot p_-(e_j)\cdot\phi.\ee
For every vector $X$ on $T\wit{M}$, we have $g\pa{p_-(X),p_-(X)}=0$ and therefore $p_-(X)\cdot p_-(X)\cdot\varphi=0$ for every $\varphi$ in $\Sigma\wit{M}_{|_M}$. Hence
\[\wih{D}\psi=-\alpha\sum_{j=1}^np_+(e_j)\cdot  p_-(e_j)\cdot\phi.\]
But, since $M$ is Lagrangian in $\wit{M}$, the complex vectors $Z_j:=p_+(e_j)$ and $\ovl{Z}_j:=p_-(e_j)$ ($1\leq j\leq n$) form a Witt-basis for $T\wit{M}\otimes\mathbb{C}$. Now remember the expression of the K\"ahler form $\wit{\Omega}$ of $(\wit{M},g,J)$ in that basis:
\[\wit{\Omega}=-2i\sum_{j=1}^nZ_j\wedge\ovl{Z}_j.\]
We deduce from that identity that
\be
\wih{D}\psi&=&-\alpha\sum_{j=1}^n\pa{Z_j\wedge\ovl{Z}_j}\cdot\phi+\alpha\sum_{j=1}^ng\pa{Z_j,\ovl{Z_j}}\phi\\
&=&-\frac{i\alpha}{2}\wit{\Omega}\cdot\phi+\frac{n\alpha}{2}\phi\\
&=&\frac{(n+1)\alpha}{2}\phi,
\ee
since $\wit{\Omega}\cdot\phi=i\phi$. A similar computation gives
\be
\wih{D}\phi&=&-\alpha\sum_{j=1}^np_-(e_j)\cdot p_+(e_j)\cdot \psi\\
&=&-\alpha\sum_{j=1}^n\pa{\ovl{Z}_j\wedge Z_j}\cdot\psi+\alpha\sum_{j=1}^ng\pa{\ovl{Z}_j,Z_j}\psi\\
&=&\frac{i\alpha}{2}\wit{\Omega}\cdot\psi+\frac{n\alpha}{2}\psi\\
&=&\frac{(n+1)\alpha}{2}\psi,\ee since $\wit{\Omega}\cdot\psi=-i\psi$. We therefore obtain: $\wih{D}^2\psi=\frac{(n+1)^2\alpha^2}{4}\psi$ and $\wih{D}^2\phi=\frac{(n+1)^2\alpha^2}{4}\phi$, i.e.,
\begin{equation}\label{ecDt}\pa{D_M^{\Sigma N}}^2\varphi=\frac{(n+1)^2\alpha^2}{4}\varphi+\frac{n^2|H|^2}{4}\varphi+\frac{n}{2}\sum_{j=1}^ne_j\cdot \nabla_{e_j}^NH\cdot\varphi\end{equation} for $\varphi:=\psi$ or $\phi$. Taking the Hermitian inner product of the sum of the two identities for $\psi$ and $\phi$ with $\psi+\phi$ and integrating lead to
\begin{equation}\label{eQRDt}\mathcal{Q}\pa{(D_M^{\Sigma N})^2,\psi+\phi}=\frac{(n+1)^2\alpha^2}{4}+\frac{n^2\int_M|H|^2\la\psi+\phi\,,\psi+\phi\ra v_g}{4\int_M\la\psi+\phi\,,\psi+\phi\ra v_g}.\end{equation}

\noindent Here we recall that, as the operator $D_M^{\Sigma N}$ is self-adjoint, we only keep the real parts when taking the Hermitian inner product of both members of (\ref{ecDt}) with $\psi+\phi$. That is why the term involving $\nabla^NH$ doesn't give any contribution to (\ref{eQRDt}).\\
To conclude, it remains to note that, if $\alpha$ is real, then the length-function of $\psi+\phi$ is constant on $\wit{M}$ (hence on $M$), whereas it cannot be constant when $\alpha$ is imaginary (see \cite{Kirch90}). From the Min-Max principle, we straightforward obtain the result.\findemo

\begin{erem}\label{rnbvalp}
{\rm If $\alpha\neq 0$, one can actually bound a greater number of eigenvalues for $D_M^{\Sigma N}$ in Theorem \ref{tmajLagr}: indeed, let $\mathcal{K}_{\frac{n-1}{2}}$ (resp. $\mathcal{K}_{\frac{n+1}{2}}$) be the (pointwise) orthogonal projection of $\mathcal{K}_\alpha$ onto the $-i$- (resp. $i$-) eigenspace of the Clif\mbox{}ford action of $\wit{\Omega}$, that is,
\be\mathcal{K}_{\frac{n-1}{2}}&:=&\left\{\psi\in\Gamma(\Sigma\wit{M})\,/\,\wit{\Omega}\cdot\psi=-i\psi\textrm{ and }\exists\, \phi\in\Gamma(\Sigma\wit{M})\,/\, (\psi,\phi)\in\mathcal{K}_\alpha\right\}\\
\mathcal{K}_{\frac{n+1}{2}}&:=&\left\{\phi\in\Gamma(\Sigma\wit{M})\,/\,\wit{\Omega}\cdot\phi=i\phi\textrm{ and }\exists\, \psi\in\Gamma(\Sigma\wit{M})\,/\, (\psi,\phi)\in\mathcal{K}_\alpha\right\}.\ee From the obvious injectivity of the orthogonal projections $\mathcal{K}_\alpha\lra\mathcal{K}_{\frac{n\pm1}{2}}$, we deduce that $\mathrm{dim}_\mathbb{C}\pa{\mathcal{K}_{\frac{n\pm1}{2}}}=\mathrm{dim}_\mathbb{C}\pa{\mathcal{K}_{\frac{n\pm1}{2}}}=N$. The identity (\ref{ecDt}) holding then on $\mathcal{K}_{\frac{n-1}{2}}\oplus\mathcal{K}_{\frac{n+1}{2}}$, it follows that, for every $\varphi$ in $\mathcal{K}_{\frac{n-1}{2}}\oplus\mathcal{K}_{\frac{n+1}{2}}$,
\[\mathcal{Q}\pa{(D_M^{\Sigma N})^2,\varphi}=\frac{(n+1)^2\alpha^2}{4}+\frac{n^2\int_M|H|^2\la\varphi\,,\varphi\ra v_g}{4\int_M\la\varphi\,,\varphi\ra v_g},\] in particular \[\mathcal{Q}\pa{(D_M^{\Sigma N})^2,\varphi}\leq\frac{(n+1)^2\alpha^2}{4}+\frac{n^2}{4}||H||_\infty^2,\] from which the following estimate holds:
\begin{equation}\label{emajnbopt}
\lambda_{2N}^2\leq\frac{(n+1)^2\alpha^2}{4}+\frac{n^2}{4}||H||_\infty^2.
\end{equation} We can therefore estimate $2N$ eigenvalues with the right member of (\ref{emajnbopt}); when $\alpha$ is a non-zero real number, we do not obtain more information on the $N$ f\mbox{}irst eigenvalues, who remain sharply bounded from Theorem \ref{tmajLagr}. However, when $\alpha$ is purely imaginary, we obtain an estimate for a larger number of eigenvalues.
}\end{erem}
\begin{ecor}\label{cmajLagr}
Under the hypotheses of Theorem \ref{tmajLagr}, assume furthermore that the complex structure $J$ induces an isomorphism between the Spin structures of $TM$ and $NM$. Then the $N$ smallest eigenvalues (counted with their multiplicities) $\lambda$ of the Hodge\,-\,de Rham Laplacian on $\Omega^{\frac{n-1}{2}}(M)\oplus \Omega^{\frac{n+1}{2}}(M)$ satisfy \[\lambda\leq\frac{(n+1)^2\alpha^2}{4}+\left\{\begin{array}{ll}\frac{n^2}{4\mathrm{Vol}(M)}\int_M|H|^2v_g&\textrm{ if }\alpha\textrm{ is real}\\ & \\ \frac{n^2||H||_\infty^2}{4}&\textrm{ if }\alpha\textrm{ is purely imaginary,}\end{array}\right.\] where $H$ is the mean curvature vector f\mbox{}ield of $M$ in $\wit{M}$. Moreover, if $M$ is minimal in $\wit{M}$, then the same result holds for the $[\frac{N+1}{2}]$ smallest eigenvalues of the Hodge\,-\,de Rham Laplacian on the space of \emph{closed} $\frac{n+1}{2}$-forms.\\
\end{ecor}

\noindent{\it Proof}: From Corollary \ref{cfS2nLL}, if $J$ identif\mbox{}ies the Spin structures of $TM$ and $NM$, then $(D_M^{\Sigma N})^2=d\delta+\delta d$. Furthermore, the isomorphism (\ref{IfS2nL}) identif\mbox{}ies the eigenspace associated to the eigenvalue $i(2p-n)$ of the Clif\mbox{}ford action of $\wit{\Omega}$ with $\Lambda^pTM\otimes\mathbb{C}$; since, under that action, the spinor f\mbox{}ield $\phi$ (resp. $\psi$) is eigen for the eigenvalue $i$ (resp. $-i$), it is a $\frac{n+1}{2}$-form (resp. a $\frac{n-1}{2}$-form) on $M$. Hence the f\mbox{}irst statement holds.\\
If moreover $H=0$, then $\wih{D}=D_M^{\Sigma N}=d+\delta$. From the equalities $\wih{D}\psi=\frac{(n+1)\alpha}{2}\phi$ and $\wih{D}\phi=\frac{(n+1)\alpha}{2}\psi$, we then deduce that
\[
\left|\begin{array}{ll}d\psi&=\frac{(n+1)\alpha}{2}\phi\\
\delta\psi&=0\end{array}\right.\qquad\textrm{ and }\qquad\left|\begin{array}{ll}d\phi&=0\\
\delta\phi&=\frac{(n+1)\alpha}{2}\psi,\end{array}\right.\] i.e., $\psi$ is coclosed and $\phi$ is closed. As the spectrum of the $\frac{n-1}{2}$-Laplacian on coclosed forms coincides with that of the $\frac{n+1}{2}$-Laplacian on closed forms (use the Hodge star operator), we obtain the second property.\findemo

\subsection{Examples}

\noindent For an odd integer $n\geq 3$, consider the round sphere $S^n$ (of constant sectional curvature 1) of dimension $n$ as canonically embedded in the $2n+1$-dimensional round sphere $S^{2n+1}$. That embedding is isometric, totally geodesic, and the canonical complex structure of $\mathbb{R}^{2n+2}$ maps the tangent bundle of $S^n$ into the horizontal space $\mathcal{H}$ def\mbox{}ined, for each $z$ in $S^{2n+1}$ as
\[\mathcal{H}_z:=\{\mathbb{R}z\oplus\mathbb{R}Jz\}^\perp\subset T_zS^{2n+1},\]
with the following property:
\[\mathcal{H}_{|_{S^n}}=TS^n\buil{\oplus}{\perp}J(TS^n).\]
Let then $\mathbb{C}\mathrm{P}^n$ be the complex projective space of complex dimension $n$. Composing the Hopf f\mbox{}ibration $S^{2n+1}\lra\mathbb{C}\mathrm{P}^n$ with the above embedding yields an immersion
\begin{equation}\label{iSnCPn}
S^n\lra\mathbb{C}\mathrm{P}^n
\end{equation}
satisfying the following: it is isometric (the Hopf f\mbox{}ibration induces an isometry from $\mathcal{H}$ onto $T\mathbb{C}\mathrm{P}^n$), totally geodesic (the Hopf f\mbox{}ibration maps horizontal geodesics onto geodesics) and Lagrangian (the Hopf f\mbox{}ibration is ``holomorphic'' w.r.t. the complex structures of $\mathcal{H}$ and $\mathbb{C}\mathrm{P}^n$). Furthermore, if $n$ is odd, the manifold $\mathbb{C}\mathrm{P}^n$ is Spin, has a unique Spin structure since it is simply-connected \cite{LM}, and carries a $2C_n^{\frac{n+1}{2}}$-dimensional space of $1$-K\"ahlerian Killing spinors \cite{Kirch90}. The round sphere $S^n$ is also Spin, and for the same reason has a unique Spin structure; more generally, every Spin vector bundle on $S^n$ has a unique Spin structure, that is, two Spin structures on a vector bundle on $S^n$ will always be isomorphic.\\
Consider then the (canonical) Spin structure of $TS^n$ and the induced one on the normal bundle of $S^n$ in $\mathbb{C}\mathrm{P}^n$ w.r.t. (\ref{iSnCPn}); then the complex structure of $\mathbb{C}\mathrm{P}^n$ will necessarily induce an isomorphism between the Spin structures of the tangent and normal bundles of $S^n$. Hence we obtain from Corollary \ref{cmajLagr} and Note \ref{rnbvalp} the existence of the following upper bound for the $2C_n^{\frac{n+1}{2}}$ smallest eigenvalues $\lambda$ of the Hodge\,-\,de Rham Laplacian on the closed $\frac{n+1}{2}$-forms:
\[\lambda\leq\frac{(n+1)^2}{4}.\]
That estimate is sharp: indeed, for $1\leq p\leq n-1$, the spectrum of the Hodge\,-\,de Rham Laplacian on the closed $p$-forms on $S^n$ is \cite{GallMey75}
\[\{(k+p)(n-p+k+1)\;/\;k\in\mathbb{N}\},\] and the multiplicity of the f\mbox{}irst eigenvalue ($k=0$) is $C_{n+1}^p$. But, for $p:=\frac{n+1}{2}$, we have $2C_n^{\frac{n+1}{2}}=C_n^{\frac{n-1}{2}}+C_n^{\frac{n+1}{2}}=C_{n+1}^{\frac{n+1}{2}}$, which is precisely the multiplicity of the f\mbox{}irst eigenvalue of the Hodge\,-\,de Rham Laplacian on the closed $\frac{n+1}{2}$-forms.\\

\noindent A further interesting example would be to consider the real $n$-dimensional projective space (with $n=4k+3$) in the complex projective space $\mathbb{C}\mathrm{P}^n$. That question, which is linked to determining \emph{all} the Lagrangian submanifolds which satisfy the equality in Theorem \ref{tmajLagr}, will be considered in a forthcoming work.

\providecommand{\bysame}{\leavevmode\hbox to3em{\hrulefill}\thinspace}

$ $\\

\noindent Institut f\"ur Mathematik\,-\,Geometrie, Universit\"at Potsdam, PF 60 15 53, 14415 Potsdam, Germany\\ ginoux@math.uni-potsdam.de

\end{document}